\documentclass{article}
\usepackage{graphicx} 

\usepackage{amsmath}
\usepackage{amssymb}
\usepackage{amsfonts}
\usepackage{amsthm}
\usepackage{mathtools}
\usepackage{enumerate}
\usepackage[hidelinks]{hyperref}
\usepackage{titlesec}
\usepackage{url}

\usepackage[T1]{fontenc}   
\usepackage{lmodern}       
\usepackage{textcomp}      

\usepackage[shortlabels]{enumitem}
\usepackage{yfonts}
\usepackage{setspace}
\usepackage[utf8]{inputenc}
\usepackage[english]{babel}
\usepackage[margin = 1in]{geometry}
\usepackage{tikz}
\usetikzlibrary{cd}
\usepackage{graphicx}
\usepackage{authblk}
\usetikzlibrary{calc}
\usepackage{fancyhdr}
\usepackage{float}

\DeclareMathOperator{\diam}{diam}
\DeclareMathOperator{\Tang}{T}

\theoremstyle{definition}

\theoremstyle{remark}
\newtheorem*{rem}{Remark}
\theoremstyle{plain}
\newtheorem{conj}{Conjecture}
\newtheorem{que}{Question}
\newtheorem{thm}{Theorem}[section]
\newtheorem{prop}[thm]{Proposition}
\newtheorem{lem}[thm]{Lemma}
\newtheorem{cor}[thm]{Corollary}
\newtheorem{defi}[thm]{Definition}

\titleformat{\section}[block]{\centering\Large\bfseries}{\thesection}{1em}{}

\title{Inscriptions in non-Euclidean Geometries}
\author{Ali Naseri Sadr}
\date{}

\begin{document}

\maketitle

\begin{abstract}
    We show how inscription problems in the plane can be generalized to Riemannian surfaces of constant curvature. We then use ideas from symplectic and Riemannian geometry to prove these generalized versions for smooth Jordan curves in the hyperbolic plane, and we prove a rectangular inscription theorem for smooth Jordan curves on the two sphere that do not intersect their antipodal. 
\end{abstract}

\section{Introduction}
Let $\gamma$ be a smooth Jordan curve in the plane. For every rectangle $R$, there are four points on $\gamma$ which are vertices of a rectangle similar to $R$; this was proved by Greene and Lobb in \cite{RectanglePaper}. 
They generalized this result to all cyclic quadrilaterals in \cite{CyclicPaper}. 
In other words, a smooth Jordan curve \textit{inscribes} every cyclic quadrilateral up to similarity.

Suppose $S$ denotes the hyperbolic plane $\mathbb{H}$ or the round sphere $S^2$. We are going to define a notion of cyclic quadrilaterals and their similarity type in $S$. Then we will show how the corresponding results in the plane can be generalized to smooth embedded closed curves in $S$ using ideas from symplectic and Riemannian geometry. 

\begin{defi}
\label{cyclic def}
    Let $(\Sigma, g)$ be an oriented Riemannian surface, and consider four points $p_1,\dots, p_4$ on $\Sigma$. We say these points are vertices of a cyclic quadrilateral with respect to the metric $g$ if they all lie on a geodesic circle. Fix three angles $\theta, \varphi_1,$ and $\varphi_2$ such that $0<\theta<\varphi_1\leq \pi$, and $\varphi_1<\varphi_2+\theta<2\pi$. The points $p_1,\dots, p_4$ make a cyclic quadrilateral of type $(\theta, \varphi_1, \varphi_2)$ if one can find $x\in\Sigma$ and $v\neq 0\in T_x\Sigma$ with
    \begin{equation*}
        p_1 = \exp_g(x,v), p_2=\exp_g(x, e^{i\theta}\cdot v), p_3=\exp_g(x, e^{i\varphi
        _1}\cdot v), p_4=\exp_g(x, e^{i(\varphi_2+\theta)}\cdot v).
    \end{equation*}
    See the following figure; if one considers the case $\varphi_1= \varphi_2=\pi$, then we call the resulting cyclic quadrilateral a rectangle of type $\theta$.

\begin{figure}[h]
        \centering
\begin{tikzpicture}[>=stealth]     
  \def\r{2}                        

  \coordinate (O) at (0,0);        
  \coordinate (P1) at ( 80:\r);    
  \coordinate (P2) at (150:\r);    
  \coordinate (P3) at (240:\r);    
  \coordinate (P4) at (290:\r);    

  \draw[blue] (O) circle (\r);
  \foreach \P/\name/\pos in {%
      P1/$p_1$/above,
      P2/$p_2$/above left,
      P3/$p_3$/below left,
      P4/$p_4$/below right}{%
        \draw[dashed, red] (O)--(\P);           
        \fill (\P) circle [radius=2pt] node[\pos] {\name};
  }

 \fill (O) circle[radius=2pt] node[below right] {$x$};

  \draw (O) ++(80:0.2) arc (80:150:0.2);
  \draw (O) ++(80: 0.6) arc (80:240:0.6);
  \draw (O) ++(150: 1.1) arc (150:290:1.1);

  \node at ($(O)+(107:0.4)$) {$\theta$};
  \node at ($(O)+(180:0.8)$) {$\varphi_1$};
  \node at ($(O)+(260:1.3)$) {$\varphi_2$};
  
\end{tikzpicture}
\caption{A cyclic quadrilateral of type $(\theta, \varphi_1, \varphi_2)$.}
\label{fig1}
\end{figure}
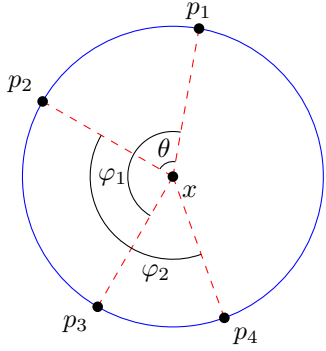
\end{defi}
In the previous definition, we used the fact that $\Sigma$ is orientable; however, if the surface is non-orientable, it is still possible to define a rectangle of type $\theta$ with respect to a Riemannian metric.
\begin{defi}
For a non-orientable Riemannian surface $(\Sigma, g)$, we say $p_1,\dots, p_4$ are four vertices of a type $\theta$ rectangle if there exist $x\in \Sigma$ and non-zero vectors $v_1, v_2\in T_x\Sigma$ with $\Vert v_1\Vert = \Vert v_2\Vert$ such that
\begin{equation*}
    p_1 = \exp_g(x, v_1), p_2=\exp(x, v_2), p_3=\exp(x, -v_1), p_4=\exp(x, -v_2),
\end{equation*}
and the angle between $v_1, v_2$ is $\theta$.
\end{defi}

\begin{defi}
    We say a closed embedded curve $\gamma$ in $(\Sigma, g)$ inscribes a cyclic quadrilateral of type $(\theta, \varphi_1, \varphi_2)$ if there are four points $p_1,\dots, p_4$ on $\gamma$ so that they are vertices of a cyclic quadrilateral of type $(\theta, \varphi_1, \varphi_2)$.
\end{defi}

\begin{thm}
\label{main hyperbolic}
    Suppose $\gamma$ is a smooth closed embedded curve in $\mathbb{H}$, and choose three angles $\theta, \varphi_1, \varphi_2$ according to Definition \ref{cyclic def}. Then $\gamma$ inscribes a cyclic quadrilateral of type $(\theta, \varphi_1, \varphi_2)$.
\end{thm}
\begin{thm}
\label{main spherical}
    Consider a smooth Jordan curve $\gamma$ in $S^2$. Let $A\colon S^2\to S^2$ denote the antipodal map, and assume $A(\gamma)\cap \gamma = \emptyset$. Then $\gamma$ inscribes a rectangle of type $\theta$ with respect to the round metric \break for every $\theta\in (0, \pi)$.
\end{thm}
\begin{cor}
\label{rect and const curvature}
    Let $(\Sigma, g)$ be a complete Riemannian surface with constant curvature other than $S^2$, and consider a smooth embedded null-homotopic curve $\gamma$ in $\Sigma$. Then $\gamma$ inscribes a rectangle of type $\theta$ for every $\theta\in(0, \pi)$.
\end{cor}
\begin{proof}
    Lift $\gamma$ to a closed curve $\tilde{\gamma}$ in the universal cover of $\Sigma$. If this universal cover is the plane, then we know $\tilde{\gamma}$ inscribes every rectangle by Greene and Lobb's result; if the universal cover is $\mathbb{H}$, then we know $\tilde{\gamma}$ inscribes every rectangle by Theorem \ref{main hyperbolic}. Finally, if the universal cover is $S^2$, then $\Sigma$ has to be $\mathbb{RP}^2$, and since $\gamma$ is an embedded curve, we must have $A(\tilde{\gamma})\cap \tilde{\gamma}=\emptyset$. Hence, we can apply Theorem \ref{main spherical} to $\tilde{\gamma}$. Now notice that a rectangle of type $\theta$ in $\tilde{\gamma}$ descends to one with the same type in $\gamma$. 
\end{proof}

If we consider the isometric embedding of $S^2$ into $\mathbb{R}^3$ as the unit sphere, then it is straightforward to check that a rectangle of type $\theta$ in $S^2$ is the same as a rectangle in $\mathbb{R}^3$ with angle $\theta$ between its two diagonals. Hence, we get the following.
\begin{cor}
    Assume $\gamma\colon S^1\to \mathbb{R}^3$ is a smooth embedding with $\Vert\gamma(t)\Vert=r$ for every $t$ and $\diam(\gamma)<2r$. Then $\gamma$ inscribes a rectangle with angle $\theta$ between its diagonals for every $\theta\in (0, \pi)$. 
\end{cor}
\begin{proof}
    This follows by scaling $\gamma$ and applying Theorem \ref{main spherical} to the unit sphere in $\mathbb{R}^3$.
\end{proof}

The outline of the paper is as follows. First, we handle the hyperbolic case and prove Theorem \ref{main hyperbolic} in Section $1$. Then we prove Theorem \ref{main spherical} in Section $2$. We conclude the paper with some remarks and questions in Section 3.
\section*{\mdseries Acknowledgments}
The author is grateful to his advisors, John Baldwin and Josh Greene, for their invaluable guidance and support. He is also thankful to Andrew Lobb and Rich Schwartz for insightful conversations on this work.

\section{Inscriptions in the Hyperbolic Plane}
\subsection{A diagonal symplectic form}
Endow the hyperbolic plane with the induced orientation and complex structure from $\mathbb{C}$; we denote this complex structure by $j$. Fix a triple of angles $(\theta, \varphi_1, \varphi_2)$ according to Definition \ref{cyclic def}. 

Let $R_\theta\colon T\mathbb{H}\to T\mathbb{H}$ denote the map
\begin{equation*}
    (x,v)\mapsto (x, e^{i\theta}\cdot v).
\end{equation*}
We define a map $F_\varphi\colon T\mathbb{H}\to\mathbb{H}\times\mathbb{H}$ by
\begin{equation*}
    (x,v)\mapsto (\exp(x,v), \exp(x, e^{i\varphi}\cdot v)).
\end{equation*}
For a pair of points $(p,q)\in \mathbb{H}\times\mathbb{H}$ and an angle $\varphi\in(0,2\pi)$, there is a unique element $(x,v)$ in $T\mathbb{H}$ such that $F_\varphi(x,v) = (p,q)$; thus $F_\varphi$ is a diffeomorphism for every $\varphi\in (0,2\pi)$.

Consider the map $I_{(\varphi_1,\varphi_2)}^\theta\colon\mathbb{H}\times\mathbb{H}\to \mathbb{H}\times\mathbb{H}$ given by $F_{\varphi_2}\circ R_\theta\circ F_{\varphi_1}^{-1}$, and let $\gamma$ be a smooth Jordan curve in $\mathbb{H}$. It follows from Definition \ref{cyclic def} that inscriptions of a cyclic quadrilateral of type $(\theta, \varphi_1,\varphi_2)$ in $\gamma$ are in one-to-one correspondence with the set
\begin{equation*}
I_{(\varphi_1,\varphi_2)}^\theta(\mathbb{T}_\gamma)\cap\mathbb{T}_\gamma\setminus \Delta_\gamma,
\end{equation*}
where $\mathbb{T}_\gamma$ denotes the torus $\gamma\times\gamma$ in $\mathbb{H}\times\mathbb{H}$, and $\Delta_\gamma$ is the diagonal curve in $\mathbb{T}_\gamma$. Note that $I_{(\varphi_1,\varphi_2)}^\theta$ is the identity on the diagonal in $\mathbb{H}\times\mathbb{H}$, so $\Delta_\gamma\subset I_{(\varphi_1,\varphi_2)}^\theta(\mathbb{T}_\gamma)\cap\mathbb{T}_\gamma$. Let $\omega$ be the hyperbolic area form on $\mathbb{H}$. There are two projection maps $\pi_1$ and $\pi_2$ from $\mathbb{H}\times\mathbb{H}\to \mathbb{H}$, and we define $\omega_i\coloneqq \pi_i^*\omega$ for $i=1,2$. Our goal is to find  non-zero real numbers $a,b,c$ such that 
\begin{equation*}
      (I_{(\varphi_1,\varphi_2)}^\theta)^*(a\omega_1+\omega_2) = b\omega_1+ c\omega_2.
\end{equation*}
If we can find such numbers, then $\mathbb{T}_\gamma$ and $I_{(\varphi_1,\varphi_2)}^\theta(\mathbb{T}_{\gamma})$ both become Lagrangians in $(\mathbb{H}\times\mathbb{H}, a\omega_1\oplus\omega_2)$, and we can reduce the inscription problem to a Lagrangian intersection problem. In order to find $a,b,c$, we first need to compute $DI_{(\varphi_1,\varphi_2)}^\theta = DF_{\varphi_2}\circ DR_\theta\circ DF_{\varphi_1}^{-1}$. 

For each element $(x,v)$ in $T\mathbb{H}$, we denote the vertical vectors in $T_{(x,v)}T\mathbb{H}$ by $\mathcal{V}_{(x,v)}$ and the horizontal ones by $\mathcal{H}_{(x,v)}$. Let $pr\colon T\mathbb{H}\to \mathbb{H}$ be the projection map; it is straightforward to see 
\begin{equation}
\label{Horiz R_t deri}
    Dpr\circ DR_\theta = Dpr
\end{equation} 
on $\mathcal{H}_{(x,v)}$.
Denote the canonical isomorphism from $\mathcal{V}_{(x,v)}$ to $T_x\mathbb{H}$ by $q_v$; then we get
\begin{equation}
\label{R deri eq}
    q_{(e^{i\theta}\cdot v)}\circ DR_\theta = R_\theta\circ q_v
\end{equation}
over the vertical vectors.

Consider a point $(x,v)\in T\mathbb{H}$ with $v\neq 0$, and suppose $(p,q) = F_\varphi(x,v)$ for some $\varphi\in (0,2\pi)$. We define two vectors
$w_p\in T_p\mathbb{H}$ and $w_q\in T_q\mathbb{H}$ by
\begin{equation}
\label{wp wq eq}
    w_p\coloneqq \frac{\frac{d}{dt}\big|_{t=1}\exp(x, t\cdot v)}{\Vert v\Vert}, \hspace{3mm} w_q\coloneqq \frac{\frac{d}{dt}\big|_{t=1}\exp(x, te^{i\varphi}\cdot v)}{\Vert v\Vert}.
\end{equation}   
 Let $\tilde{v} = \frac{v}{\Vert v\Vert}$, and denote the horizontal lifts of $\tilde{v}$ and $j\tilde{v}$ in $\mathcal{H}_{(x,v)}$ by $\tilde{v}_H$ and $j\tilde{v}_H$; the vectors $\langle \tilde{v}, j\tilde{v}, \tilde{v}_H, j\tilde{v}_H\rangle$ give a basis for $T_{(x,v)}T\mathbb{H}$ where we have identified the vertical vectors with $T_x\mathbb{H}$.

We will need the following lemma for later computations.
\begin{lem}
\label{Exp Deri Compute}
    Let $(x,v)$ be an element in $T\mathbb{H}$ with $p = \exp(x,v)$. Define $w_p$ as above; then we have
    \begin{align*}
        & D\exp_{(x,v)}(\tilde{v}_H) = w_p,\\
         D\exp&_{(x,v)}(j\tilde{v}_H) = \cosh(\Vert v\Vert)\cdot jw_p.
    \end{align*}
\end{lem}
\begin{proof}
    The first equality follows from the definition of $w_p$. For the second one, consider a curve $x(s)$ with $x(0) = x$ and $\Dot{x}(0) = j\tilde{v}$. Denote the parallel transport of $v$ along $x(s)$ by $v(s)$. Consider $\Gamma\colon\mathbb{R}^2\to \mathbb{H}$ defined by $\Gamma(s,t) = \exp(x(s), t\cdot v(s))$. Let $u(t) = \frac{\partial\Gamma}{\partial s}(0, t)$;  we have
    $D\exp_{(x,v)}(j\tilde{v}_H) = u(1)$. Note that $u(t)$ is a Jacobi vector field along the geodesic $\exp(x, t\cdot v)$, and 
    \begin{equation*}
        \frac{Du}{Dt}(0) = \frac{D}{Dt}\frac{\partial\Gamma}{\partial s}(0, 0)=\frac{D}{Ds}\frac{\partial\Gamma}{\partial t}(0,0)=\frac{Dv}{\partial s}(0) = 0.
    \end{equation*}
    It follows from the Jacobi equation for $\mathbb{H}$ that $u(t) = \cosh(t\Vert v\Vert)\cdot jr(t)$ where $r(t)$ is the unit tangent vector to $\exp(x, t\cdot v)$; we refer the unfamiliar reader to \cite{RiemBook}.
    Hence, we get $u(1) = \cosh(\Vert v\Vert)\cdot jw_p$.
\end{proof}
\begin{lem}
\label{DF_phi Computation}
    The map $DF_\varphi\colon T_{(x,v)}T\mathbb{H}\to T_p\mathbb{H}\oplus T_q\mathbb{H}$ is given by
    \begin{align}
         &\tilde{v} \mapsto (w_p, w_q), \label{F_varphi Deri 1}\\
          j\tilde{v}\mapsto &\frac{\sinh(\Vert v\Vert)}{\Vert v\Vert}\cdot(jw_p, jw_q)\label{F_varphi Deri 2},
    \end{align}
    and
    \begin{align}
        & \tilde{v}_H\mapsto (w_p, \cos(\varphi)\cdot w_q-\sin(\varphi)\cosh(\Vert v\Vert)\cdot jw_q),\label{F Deri 3} \\
        j\tilde{v}_H& \mapsto (\cosh(\Vert v\Vert)\cdot jw_p, \sin(\varphi)\cdot w_q+\cos(\varphi)\cosh(\Vert v\Vert)\cdot jw_q)\label{F Deri 4}.
    \end{align}
\end{lem}
\begin{proof}
    The first equality follows from the definition of $w_p$ and $w_q$. For the second one, we have
    \begin{equation*}
        DF_\varphi(j\tilde{v}) = (u_p, u_q) = \frac{d}{ds}\Big|_{s=0}\Big(\exp(x, e^{is}\cdot v), \exp(x, e^{i(s+\varphi)}\cdot v)\Big),
    \end{equation*}
    where $u_p\in T_p\mathbb{H}$ and $u_q\in T_q\mathbb{H}$. Both $u_p$ and $u_q$ are tangent to the geodesic circle with radius $\Vert v\Vert$ around $x$; hence they are orthogonal to $w_p$ and $w_q$ respectively. Moreover, $\Vert u_p\Vert=\Vert u_q\Vert$ because there is an isometry of $\mathbb{H}$ that takes one to the other. Therefore, we must have $u_p = c\cdot jw_p$ and $u_q = c\cdot jw_q$; the constant $c$ is equal to $\frac{\sinh\Vert v\Vert}{\Vert v\Vert}$ because the geodesic circle have length $2\pi\sinh(\Vert v\Vert)$. We have
    \begin{equation}
    \label{F_p horiz eq}
        DF_\varphi(\tilde{v}_H) = (D\exp_{(x,v)}(\tilde{v}_H), D\exp\circ DR_\varphi(\tilde{v}_H)). 
    \end{equation}
    The first component is equal to $w_p$. We denote $e^{i\varphi}\cdot v$ by $z$, and we have 
    \begin{equation*}
        DR_\varphi(\tilde{v}_H) = \cos(\varphi)\cdot\tilde{z}_H -\sin(\varphi)\cdot j\tilde{z}_H.
    \end{equation*}
    Now it follows from Lemma \ref{Exp Deri Compute} that the second component in equation \eqref{F_p horiz eq} is given by 
    \begin{equation*}
        \cos(\varphi)\cdot w_q-\sin(\varphi)\cosh(\Vert v\Vert)\cdot jw_q
    \end{equation*}
    Similarly, we get
    \begin{equation}
    \label{F_p horiz eq 2}
        DF_\varphi(j\tilde{v}_H) = (D\exp_{(x,v)}(j\tilde{v}_H), D\exp\circ DR_\varphi(j\tilde{v}_H)). 
    \end{equation}
    The first component is equal to $\cosh(\Vert v\Vert)\cdot jw_p$ by the previous lemma; for the second component, we have
    \begin{equation*}
         DR_\varphi(j\tilde{v}_H) = \sin(\varphi)\cdot\tilde{z}_H +\cos(\varphi)\cdot j\tilde{z}_H.
    \end{equation*}
    Therefore, the second component in \eqref{F_p horiz eq 2} is equal to
    \begin{equation*}
        \sin(\varphi)\cdot w_q+\cos(\varphi)\cosh(\Vert v\Vert)\cdot jw_q
    \end{equation*}
    by Lemma \ref{Exp Deri Compute}.
\end{proof}

Let $(p, q)\in \mathbb{H}\times\mathbb{H}$ be a pair of points with $p\neq q$, and suppose $(p', q') = I_{(\varphi_1, \varphi_2)}^\theta(p,q)$. Choose $(x,v)\in T\mathbb{H}$ such that
\begin{align*}
    & (p, q)=(\exp(x,v), \exp(x, e^{i\varphi_1}\cdot v)), \\
     (p', &q')=(\exp(x, e^{i\theta}\cdot v), \exp(x, e^{i(\varphi_2+\theta)}\cdot v)).
\end{align*}
Define $w_p$ and $w_q$ by equation \eqref{wp wq eq}, and let
\begin{equation*}
     w_{p'}\coloneqq \frac{\frac{d}{dt}\big|_{t=1}\exp(x, te^{i\theta}\cdot v)}{\Vert v\Vert}, \hspace{3mm} w_{q'}\coloneqq \frac{\frac{d}{dt}\big|_{t=1}\exp(x, te^{i(\varphi_2+\theta)}\cdot v)}{\Vert v\Vert}.
\end{equation*}
\begin{lem}
\label{DI Comp Lemma}
    The map $DI_{(\varphi_1, \varphi_2)}^{\theta}\colon T_p\mathbb{H}\oplus T_q\mathbb{H}\to T_{p'}\mathbb{H}\oplus T_{q'}\mathbb{H}$ is given by
    \begin{equation}
    \label{I deri 1}
        (w_p, w_q) \mapsto (w_{p'}, w_{q'}),\hspace{3mm} (jw_p, jw_q) \mapsto (jw_{p'}, jw_{q'}).
    \end{equation}
Moreover, we have
\begin{equation}
\label{I Deri 2}
\begin{split}
        &(w_p, \cos(\varphi_1)\cdot w_q-\sin(\varphi_1)\cosh(\Vert v\Vert)\cdot jw_q)\mapsto\\ (\cos(\theta)\cdot w_{p'}&-\sin(\theta)\cosh(\Vert v\Vert)\cdot jw_{p'}, \cos(\theta+\varphi_2)w_{q'}-\sin(\theta+\varphi_2)\cosh(\Vert v\Vert)\cdot jw_{q'}), 
\end{split}
\end{equation}
and
\begin{equation}
\label{I Deri 3}
\begin{split}
     & (\cosh(\Vert v\Vert)\cdot jw_p, \sin(\varphi_1)\cdot w_q+\cos(\varphi_1)\cosh(\Vert v\Vert)\cdot jw_q)\mapsto\\
         (\sin(\theta)\cdot w_{p'}&+\cos(\theta)\cosh(\Vert v\Vert)\cdot jw_{q'}, \sin(\theta+\varphi_2)\cdot w_{q'}+\cos(\theta+\varphi_2)\cosh(\Vert v\Vert)\cdot jw_{q'}).
\end{split}
\end{equation}
\end{lem}
\begin{proof}
We know $DI_{(\varphi_1, \varphi_2)}^{\theta} = DF_{\varphi_2}\circ DR_\theta\circ DF_{\varphi_1}^{-1}$, so equation \eqref{I deri 1} follows from equations \eqref{F_varphi Deri 1} and \eqref{F_varphi Deri 2} combined with equation \eqref{R deri eq}. For equation \eqref{I Deri 2}, note that we have
\begin{equation*}
    DF_{\varphi_1}^{-1}((w_p, \cos(\varphi_1)\cdot w_q-\sin(\varphi_1)\cosh(\Vert v\Vert)\cdot jw_q)) = \tilde{v}_H
\end{equation*}
by Lemma \ref{DF_phi Computation}. Let $\tilde{v}_H^{\theta}$ be the horizontal lift of $e^{i\theta}\cdot \tilde{v}$ in $\mathcal{H}_{(x, e^{i\theta}\cdot v)}$ and $j\tilde{v}_H^{\theta}$ the horizontal lift of $e^{i\theta}\cdot j\tilde{v}$; we get
\begin{equation*}
    DR_\theta(\tilde{v}_H) = \cos(\theta)\cdot\tilde{v}_H^{\theta} - \sin(\theta)\cdot j\tilde{v}_H^{\theta}
\end{equation*}
by equation \eqref{Horiz R_t deri}. Therefore, we have
\begin{equation*}
    DF_{\varphi_2}\circ DR_\theta(\tilde{v}_H) = \cos(\theta)\cdot DF_{\varphi_2}(\tilde{v}_H^{\theta})-\sin(\theta)\cdot DF_2(j\tilde{v}_H^{\theta}).
\end{equation*}
Now we can deduce \eqref{I Deri 2} by applying Lemma \ref{DF_phi Computation} to the right hand side of the above identity. Equation \eqref{I Deri 3} follows in a similar way.
\end{proof}
Let $r_\varphi$ denote the following $2\times2$ matrix
\begin{equation*}
\begin{bmatrix}
    \cos(\varphi) & \sin(\varphi) \\
    -\sin(\varphi) & \cos(\varphi)
\end{bmatrix}.
\end{equation*}
We define two $4\times4$ matrices $A_{\varphi_1}$ and $M_{(\varphi_2, \theta)}$ as follows.
\begin{equation*}
A_{\varphi_1} =
    \begin{bmatrix}
        I & I \\
        I & r_{\varphi_1}
    \end{bmatrix}, \hspace{3mm}
    M_{(\varphi_2, \theta)} = 
    \begin{bmatrix}
        I & r_{\theta} \\
        I & r_{(\theta+\varphi_2)}
    \end{bmatrix}.
\end{equation*}
The vectors $\langle w_p, \cosh(\Vert v\Vert)\cdot jw_p, w_q, \cosh(\Vert v\Vert)\cdot jw_q\rangle$ make a basis for $T_p\mathbb{H}\oplus T_{q}\mathbb{H}$, and we can define a basis for $T_{p'}\mathbb{H}\oplus T_{q'}\mathbb{H}$ similarly. The map $DI_{(\varphi_1, \varphi_2)}^{\theta}$ is determined by the matrix
\begin{equation}
\label{derivative to matrix lemma}
    N_{(\theta, \varphi_1,\varphi_2)}= M_{(\varphi_2, \theta)}\cdot A_{\varphi_1}^{-1}
\end{equation}
in these bases according to Lemma \ref{DI Comp Lemma}.
\begin{lem}
\label{matrix pullback lemma}
    There exist non-zero constants $a(\theta, \varphi_1,\varphi_2), b(\theta, \varphi_1,\varphi_2),$ and $c(\theta, \varphi_1,\varphi_2)$ such that
    \begin{equation}
    \label{pullback eq}
        N_{(\theta, \varphi_1,\varphi_2)}^*(adx_1\wedge dy_1+dx_2\wedge dy_2) = bdx_1\wedge dy_1+cdx_2\wedge dy_2.
    \end{equation}
\end{lem}
\begin{proof}
We let 
\begin{equation*}
    a = \frac{\sin(\frac{\varphi_2+\theta}{2})\sin(\frac{\varphi_2+\theta-\varphi_1}{2})}{\sin(\frac{\theta}{2})\sin(\frac{\varphi_1-\theta}{2})}.
\end{equation*}
It is easy to verify $a\neq 0$ whenever the triple $(\theta, \varphi_1, \varphi_2)$ satisfies the conditions in Definition \ref{cyclic def}. Consider the matrix $J=\begin{bmatrix}
    0 & 1\\
    -1 & 0
\end{bmatrix}$, and define
\begin{equation*}
    J_a=
    \begin{bmatrix}
        aJ & 0 \\
        0 & J
    \end{bmatrix}.
\end{equation*}
The form on the left hand side of equation \eqref{pullback eq} is given by
\begin{equation*}
    (w_1, w_2)\mapsto \langle w_1, N_{(\theta, \varphi_1,\varphi_2)}^{T}J_aN_{(\theta, \varphi_1,\varphi_2)}(w_2)\rangle.
\end{equation*}
We used the Python package SymPy to verify the matrix $K = N_{(\theta, \varphi_1,\varphi_2)}^{T}J_aN_{(\theta, \varphi_1,\varphi_2)}$ has the form 
\begin{equation*}
    \begin{bmatrix}
        bJ & 0 \\
        0 & cJ
    \end{bmatrix}
\end{equation*}
for some real numbers $b,c$ depending on $(\theta,\varphi_1,\varphi_2)$, and this proves the equation; the Python code can be found in \cite{code}.
\end{proof}
\begin{prop}
\label{pullback map prop and Lags}
    Let $\omega_a= a\omega_1+\omega_2$ be a symplectic form on $\mathbb{H}\times\mathbb{H}$ where $a$ is determined by Lemma \ref{matrix pullback lemma}. Then we have
    \begin{equation*}
        (I_{(\varphi_1, \varphi_2)}^{\theta})^* \omega_a = b\omega_1+c\omega_2.
    \end{equation*}
    In particular, if $\gamma$ is a smooth Jordan curve in $\mathbb{H}$, then both $\mathbb{T}_\gamma$ and $I_{(\varphi_1, \varphi_2)}^{\theta}(\mathbb{T}_\gamma)$ are Lagrangians in \break $(\mathbb{H}\times\mathbb{H}, \omega_a)$.
\end{prop}
\begin{proof}
    Consider a pair of points $(p,q)\in \mathbb{H}\times\mathbb{H}$, and assume $(p',q') = I_{(\varphi_1, \varphi_2)}^{\theta}(p,q)$. We work with the basis 
    $\langle w_p, \cosh(\Vert v\Vert)\cdot jw_p, w_q, \cosh(\Vert v\Vert)\cdot jw_q\rangle$ for $T_p\mathbb{H}\oplus T_q\mathbb{H}$ and the corresponding one for $T_{p'}\mathbb{H}\oplus T_{q'}\mathbb{H}$. The form $\omega_a$ is equal to
    \begin{equation*}
        \cosh(\Vert v\Vert)\cdot(adx_1\wedge dy_1+dx_2\wedge dy_2)
    \end{equation*}
    with respect to both of these bases. Hence, we get
    \begin{equation*}
        (I_{(\varphi_1, \varphi_2)}^{\theta})^* \omega_a = \cosh(\Vert v\Vert)\cdot(bdx_1\wedge dy_1+cdx_2\wedge dy_2) = b\omega_1+c\omega_2
    \end{equation*}
    by Lemma \ref{matrix pullback lemma}. For the second part, notice that $\mathbb{T}_\gamma$ is a Lagrangian with respect to every linear combination of $\omega_1$ and $\omega_2$.  
\end{proof}
\subsection{Lagrangian smoothing and Maslov number}
Consider a smooth Jordan curve $\gamma$ in $\mathbb{H}$ and a fixed triple of angles $(\theta, \varphi_1, \varphi_2)$. We fix the symplectic form $\omega_a$ given in Proposition \ref{pullback map prop and Lags}, and denote the two Lagrangian tori $\mathbb{T}_\gamma$ and $I_{(\varphi_1, \varphi_2)}^{\theta}(\mathbb{T}_\gamma)$ by $T_1$ and $T_2$ respectively.
\begin{lem}
\label{clean int lemma}
    The two Lagrangians $T_1$ and $T_2$ intersect cleanly along the diagonal loop $\Delta_\gamma$.
\end{lem}
\begin{proof}
    Fix a point $p$ on the curve $\gamma$; it is straightforward to check $DF_{\varphi_1}$ takes a vertical vector $v\in\mathcal{V}_{(p,0)}$ to $(v, e^{i\varphi_1}\cdot v)$, and it maps a Horizontal vector $v_H$ to $(v,v)$. Similarly, $DF_{\varphi_2}\circ DR_\theta$ takes a vertical vector $v$ to $(e^{i\theta}\cdot v, e^{i(\theta+\varphi_2)}\cdot v)$ and a horizontal vector $v_H$ to $(v, v)$. Using this, we get that $DI_{(\varphi_1, \varphi_2)}^{\theta}\colon T_p\mathbb{H}\oplus T_p\mathbb{H}\to T_p\mathbb{H}\oplus T_p\mathbb{H}$ acts by 
    \begin{equation*}
        (w, w) \mapsto (w,w),\hspace{3mm} (w, -w)\mapsto \Big(\frac{1+e^{i\varphi_1}-2e^{i\theta}}{e^{i\varphi_1}-1}\cdot v, \frac{1+e^{i\varphi_1}-2e^{i(\theta+\varphi_2)}}{e^{i\varphi_1}-1}\cdot v\Big)
    \end{equation*}
    for every $w\in \Tang_p\mathbb{H}$. Let $u$ be a vector tangent to $\gamma$ at $p$; then $T_{(p,p)}T_1$ is the vector space $\mathbb{R}u\oplus\mathbb{R}u$, and $T_{(p,p)}T_2$ is the image of this vector space under $DI_{(\varphi_1, \varphi_2)}^{\theta}$. It follows form the previous computation that 
    \begin{equation*}
        T_{(p,p)}T_1\cap T_{(p,p)}T_2 = \mathbb{R}(u, u)= T_{(p,p)}\Delta_{\gamma}.
    \end{equation*}
\end{proof}
Since the intersection of $T_1$ and $T_2$ is clean along $\Delta_\gamma$, we can perform a Lagrangian smoothing along the diagonal curve; this results in an immersed  Lagrangian torus $T$ where the set of self intersections of $T$ is in one-to-one correspondence with the set
\begin{equation*}
    T_1\cap T_2\setminus\Delta_\gamma.
\end{equation*}
See \cite[p.3]{CyclicPaper}
for details of this operation.
\begin{lem}
    The Lagrangian torus $T$ has Maslov number $4$.
\end{lem}
\begin{proof}
    This can be proved using an argument identical to the one in \cite{CyclicPaper}.
\end{proof}
\begin{proof}[Proof of Theorem \ref{main hyperbolic}]
    The theorem is equivalent to proving the set of self intersections of $T$ is non-empty. By contradiction, suppose this is not the case, and $T$ is an embedded Lagrangian torus in $(\mathbb{H}\times\mathbb{H}, \omega_a)$. Notice that $\mathbb{H}$ with any multiple of the hyperbolic area form is symplectomorphic to $(\mathbb{R}^2, dx\wedge dy)$. Therefore, $(\mathbb{H}\times\mathbb{H}, \omega_a)$ is symplectomorphic to $(\mathbb{C}^2, \omega_{std})$, but this gives a contradiction since every embedded Lagrangian torus in $(\mathbb{C}^2, \omega_{std})$ has Maslov number $2$; see \cite{PoltPaper, VitPaper}
    for more details.
\end{proof}

\section{Inscriptions in the Round Sphere}
\subsection{Rectangles as a Hamiltonian Motion}
Let $\omega$ denote the Riemannian area form on the round two sphere; this is a multiple of Fubini-Study form on $S^2$. Endow $S^2\times S^2$ with the symplectic form $\omega\oplus\omega$. 

We denote the distance with respect to the round metric by $d\colon S^2\times S^2\to [0,\pi]$ and the antipodal map on $S^2$ by $A$. Consider the set of antipodal pairs $D_A$ in $S^2\times S^2$; the function $d^2$ is smooth on $S^2\times S^2\setminus D_A$. If $(p,q)$ is a pair of points in the complement of $D_A$, then we can find a unique element $(x,v)\in TS^2$ with $\Vert v\Vert<\frac{\pi}{2}$, $p = \exp(x,v)$, and $q=\exp(x, -v)$. One can see that gradient of $\frac{d^2}{2}$ at $(p,q)$ with respect to the product metric on $S^2\times S^2$ is given by 
\begin{equation}
\label{wp wq eq}
    (w_p, w_q)\coloneqq\frac{d}{dt}\Big|_{t=1}\Big(\exp(tv), \exp(-tv)\Big).
\end{equation}
Define a Hamiltonian $H\colon S^2\times S^2\setminus D_A\to \mathbb{R}$ by $H(z,w)=-4\cos\big(\frac{d(z,w)}{2}\big)$. This is a smooth function because $d^2$ is smooth.
\begin{lem}
\label{Rect Motion Lem}
    Fix an angle $\theta\in(0, \pi)$, and denote the time $\theta$ Hamiltonian flow of $H$ by $\psi_{\theta}$. Assume $(p_1, p_3)$ and $(p_2, p_4)$ are in $S^2\times S^2\setminus (D_A\cup \Delta_{S^2})$.
    Then the points $p_1,p_2,p_3$, and $p_4$ are vertices of a rectangle with type $\theta$ if and only if $(p_2, p_4)=\psi_{\theta}(p_1, p_3)$.
\end{lem}
\begin{proof}
    For a distinct pair of points $(p, q)$, we have
    \begin{equation*}
        \nabla H(p,q) = 2\sin\Big(\frac{d(p,q)}{2}\Big)\cdot \frac{\nabla(\frac{d^2}{2})}{d(p,q)} = \frac{\sin(\Vert v\Vert)}{\Vert v\Vert}\cdot(w_p, w_q)
    \end{equation*}
    where $(w_p, w_q)$ and $(x, v)$ are defined in equation \eqref{wp wq eq}. Now let $j$ be the complex structure on $S^2$, and denote the complex structure $j\oplus j$ on $S^2\times S^2$ by $J$.
    Hence, we get
    \begin{equation*}
        X_H(p,q)= J\nabla H(p,q) = \frac{\sin(\Vert v\Vert)}{\Vert v\Vert}\cdot(jw_p, jw_q)
    \end{equation*}
    where $X_H$ is the Hamiltonian vector field of $H$. Thus both projections of $X_H(p,q)$ in $T_pS^2$ and $T_qS^2$ are tangent to the geodesic circle centered around $x$ and going through $p$ and $q$. The geodesic circle has length equal to $2\pi\sin(\Vert v\Vert)$, so we conclude
    \begin{equation*}
        \psi_{\theta}(p, q) = (\exp(x, e^{i\theta}\cdot v), \exp(x, -e^{i\theta}\cdot v)).
    \end{equation*}
    Note for a diagonal point $(p, p)$ in $S^2\times S^2$, we have $\psi_{\theta}(p,p) = (p,p)$ for every $\theta$.
\end{proof}

\begin{rem}
    If we consider the round sphere as the unit sphere in $\mathbb{R}^3$, then we get
    \begin{equation*}
        H(z,w)=-4\sqrt{1-\frac{\Vert z-w\Vert^2}{4}}.
    \end{equation*}
\end{rem}

\subsection{A Lagrangian Klein Bottle in $\mathbb{CP}^2$}
Consider the map $\sigma\colon S^2\times S^2\to S^2\times S^2$ that interchanges the two coordinates; this map is a symplectomorphism, and we have
\begin{equation*}
    H = H\circ\sigma.
\end{equation*}
Therefore, $\psi_\theta\circ\sigma = \sigma\circ\psi_\theta$. The quotient of $S^2\times S^2$ by $\sigma$ is $\mathbb{CP}^2$, and we denote the quotient map by $\pi$. Let $\gamma$ be a smooth Jordan curve in $S^2$ with $A(\gamma)\cap \gamma\neq \emptyset$. Then $\mathbb{T}_\gamma = \gamma\times\gamma$ and $\psi_\theta(\mathbb{T}_\gamma)$ are Lagrangians, and they are invariant under $\sigma$. The following can be proved by a computation analogous to Lemma \ref{clean int lemma}.  
\begin{lem}
    The two tori $\mathbb{T}_\gamma$ and $\psi_{\theta}(\mathbb{T}_\gamma)$ intersect cleanly along the diagonal curve $\Delta_\gamma$.
\end{lem}
Since the tori are invariant under $\sigma$, we can perform a Lagrangian smoothing of $\mathbb{T}_\gamma$ and $\psi_\theta(\mathbb{T}_\gamma)$ in a neighborhood of $\Delta_\gamma$ so that the resulting immersed Lagrangian torus $T$ is also invariant under $\sigma$, and it does not intersect the diagonal in $S^2\times S^2$. For details of this operation, we refer the unfamiliar reader to \cite{RectanglePaper}.\break
The self intersection set of $T$ is in one-to-one correspondence with the the inscriptions of type $\theta$ rectangles in $\gamma$ by Lemma \ref{Rect Motion Lem}.
\begin{proof}[Proof of Theorem \ref{main spherical}]
   We are going to prove the self intersection set of $T$ is non-empty. By contradiction, suppose $T$ is an embedded Lagrangian torus in $(S^2\times S^2, \omega\oplus\omega)$, and let $N$ be a neighborhood of the diagonal in $S^2\times S^2$ such that $T\cap N=\emptyset$. By Theorem $1$ in \cite{perutz2008remarkkahlerformssymmetric},
   we can find a symplectic form $\omega'$ on $\mathbb{CP}^2$ so that
   \begin{equation*}
       \pi^*\omega' = \omega\oplus\omega
   \end{equation*}
   in $S^2\times S^2\setminus N$. Hence, $\pi(T)$ becomes a Lagrangian Klein bottle in $(\mathbb{CP}^2, \omega')$ which gives us a contradiction because there are no embedded Lagrangian Klein bottle in a symplectic $\mathbb{CP}^2$; this was proved \break in Theorem $2$ of \cite{LagKleinCp}; see also \cite{Nemirovski_Paper}.
\end{proof}

\section{Conclusion}
We conclude the paper by mentioning some questions and conjectures for future directions. 

One can also realize the rectangles in $\mathbb{H}$ as a Hamiltonian motion in $(\mathbb{H}\oplus\mathbb{H}, \omega\oplus\omega)$; let $d$ denote the hyperbolic distance function and define
\begin{equation*}
    H(z,w) = 4\cosh\Big(\frac{d(z,w)}{2}\Big).
\end{equation*}
It is straightforward to show the Hamiltonian flow of $H$ at time $\theta\in (0, \pi)$ realizes all the rectangles of type $\theta$ in $\mathbb{H}$ analogously to Lemma \ref{Rect Motion Lem}. We believe it is possible to use this Hamiltonian motion or the one on $S^2$ in combination with the ideas of \cite{Jordan-Floer} to prove rectangular inscriptions for rectifiable Jordan curves in $\mathbb{H}$ and $S^2$.
For proving Theorem \ref{main spherical}, we had to restrict to Jordan curves in $S^2$ that are disjoint from their antipodals, since the Hamiltonian motion in Lemma \ref{Rect Motion Lem} only generates the rectangles that are not contained on a great circle. A Jordan curve intersects its antipodal if and only if it has diameter $\pi$ with respect to the round metric; we had the idea to approximate such a Jordan curve with Jordan curves that have diameter less than $\pi$, and deduce the result from Theorem \ref{main spherical}. However, it is not possible to approximate every Jordan curve of diameter $\pi$ in $S^2$ with Jordan curves that have diameter less than $\pi$. We suspect Theorem \ref{main spherical} is true for Jordan curves of diameter $\pi$; another direction for proving this result is to use ideas of Floer homology applied to the complement of the antipodal pairs on the curve in the product torus inside $S^2\times S^2\setminus D_A$. In particular, we propose the following conjecture.
\begin{conj}
    Suppose $\gamma$ is a smooth Jordan curve in $S^2$ with diameter $\pi$. Then $\gamma$ inscribes a rectangle of type $\theta$ for every $\theta\in (0, \pi)$.
\end{conj}
We were able to prove our main theorems because of the interaction between the Riemannian geometry of a surface $S$ with constant curvature and the symplectic geometry of diagonal forms on $S\times S$. It would be interesting to find out if similar ideas can be applied to other Riemannian surfaces. In particular, we think complete surfaces with non-positive curvature can be an interesting future direction.
\begin{que}
    Is it possible to prove a rectangular or cyclic inscription theorem for null-homotopic Jordan curves in complete surfaces with non-positive curvature?
\end{que}
The result in Corollary \ref{rect and const curvature} is only true for null-homotopic curves. For instance, the image of a half great circle in $\mathbb{RP}^2$ does not inscribe any rectangle; similarly, one can find essential curves in hyperbolic surfaces that do not inscribe any rectangle. However, it might be possible to prove an inscription theorem for a pair of disjoint essential Jordan curves $\gamma_1$ and $\gamma_2$ in a constantly curved surface if their lifts to the universal cover respect some symmetries depending on their homotopy classes. In particular, we have the periodic inscriptions in $\mathbb{R}^2$ (see \cite{sadr2025periodicinscriptionisoscelestrapezoids}), and we can ask the following.
\begin{que}
    What is the analog of periodic inscriptions in $\mathbb{H}$?
\end{que}
It is possible to prove a result similar to Proposition \ref{pullback map prop and Lags} for the round metric on the two sphere; however, in order to prove a cyclic inscription theorem for the two sphere, one needs to prove every Lagrangian torus in a non-monotone symplectic $S^2\times S^2$ has Maslov number two; compare to \cite{EvansPaper, Maslov2Number}.
\begin{que}
    Let $a$ be a non-zero real number. Does every Lagrangian torus in $(S^2\times S^2, a\omega\oplus\omega)$ have Maslov number $2$?
\end{que}

\bibliographystyle{plain}
\bibliography{ref}

\begin{thebibliography}{10}

\bibitem{EvansPaper}
Nikolas Adaloglou and Jonathan~David Evans.
\newblock Lagrangian klein bottles in $s^2 \times s^2$, 2024.

\bibitem{Maslov2Number}
Mihai Damian.
\newblock Floer homology on the universal cover, {A}udin's conjecture and other
  constraints on {L}agrangian submanifolds.
\newblock {\em Comment. Math. Helv.}, 87, 2012.

\bibitem{RiemBook}
Manfredo Perdig\~ao do~Carmo.
\newblock {\em Riemannian geometry}.
\newblock Mathematics: Theory \& Applications. Birkh\"auser Boston, Inc.,
  Boston, MA, 1992.

\bibitem{RectanglePaper}
Joshua~Evan Greene and Andrew Lobb.
\newblock The rectangular peg problem.
\newblock {\em Ann. of Math. (2)}, 194:509--517, 2021.

\bibitem{CyclicPaper}
Joshua~Evan Greene and Andrew Lobb.
\newblock Cyclic quadrilaterals and smooth {J}ordan curves.
\newblock {\em Invent. Math.}, 234:931--935, 2023.

\bibitem{Jordan-Floer}
Joshua~Evan Greene and Andrew Lobb.
\newblock Floer homology and square pegs, 2024.

\bibitem{Nemirovski_Paper}
S.~Yu. Nemirovski\u~i.
\newblock The homology class of a {L}agrangian {K}lein bottle.
\newblock {\em Izv. Ross. Akad. Nauk Ser. Mat.}, 73, 2009.

\bibitem{perutz2008remarkkahlerformssymmetric}
Tim Perutz.
\newblock A remark on k\"ahler forms on symmetric products of riemann surfaces,
  2008.

\bibitem{PoltPaper}
L.~V. Polterovich.
\newblock The {M}aslov class of the {L}agrange surfaces and {G}romov's
  pseudo-holomorphic curves.
\newblock {\em Trans. Amer. Math. Soc.}, 325:241--248, 1991.

\bibitem{sadr2025periodicinscriptionisoscelestrapezoids}
Ali~Naseri Sadr.
\newblock Periodic inscription of isosceles trapezoids, 2025.

\bibitem{code}
Ali~Naseri Sadr.
\newblock Sympy verification of lemma 2.4.
\newblock
  \url{https://github.com/alins95/Computation-of-the-Constant-a-in-the-paper-Inscriptions-in-non-Euclidean-Geometries-.git},
  2025.

\bibitem{LagKleinCp}
V.~V. Shevchishin.
\newblock Lagrangian embeddings of the {K}lein bottle and the combinatorial
  properties of mapping class groups.
\newblock {\em Izv. Ross. Akad. Nauk Ser. Mat.}, 73:153--224, 2009.

\bibitem{VitPaper}
C.~Viterbo.
\newblock A new obstruction to embedding {L}agrangian tori.
\newblock {\em Invent. Math.}, 100:301--320, 1990.

\end{thebibliography}

\begin{flushleft}
    Boston College. Massachusetts, USA.\\
    naserisa@bc.edu 
\end{flushleft}
\end{document}